\documentclass{article}

\usepackage{amssymb,dsfont,latexsym}

\newtheorem{thm}{Theorem}
\newtheorem{defin}[thm]{Definition}
\newtheorem{prop}[thm]{Proposition}
\newtheorem{lem}[thm]{Lemma}
\newtheorem{cor}[thm]{Corollary}

\newcommand\enu[1]{\smallskip\newline\makebox[5mm][l]{\rm(#1)}}

\newcommand\bp{\noindent{\it Proof.}\ }

\newcommand\1[1]{{\cal #1}}

\newcommand\2{{\frac{1}{2}}}

\begin{document}

\title{Multiplicative properties of positive maps}

\author{Erling St{\o}rmer}

\date{ }

\maketitle

\begin{abstract}
Let $\phi$ be a positive unital normal map of a von Neumann
algebra $M$ into itself. It is shown that there exists a largest
Jordan subalgebra $C_\phi$ of $M$ such that the restriction of
$\phi$ to $C_\phi$ is a Jordan automorphism and each weak limit
point of $(\phi^n (a))$ for $a\in M$ belongs to $C_\phi$.
\end{abstract}

\centerline{\small Dedicated to the memory of Gert K. Pedersen}

\section*{1. Introduction}

In the study of positive linear maps of $C^*$-algebras the
multiplicative properties of such maps have been studied by
several authors, see
e.g.\cite{K2},\cite{B},\cite{C1},\cite{C2},\cite{ER}. If
$\phi\colon A \to B$ is a positive unital map between
$C^*$-algebras $A$ and $B$ an application of Kadison's Schwarz
inequality,\cite{K1} to the operators $a+a^*$ and $i(a-a^*)$
yields the inequality \cite{S}
\begin{equation}\label{e1}
\phi(a\circ a^*)\geq \phi(a)\circ\phi(a)^*, a\in A,
\end{equation}
where $a\circ b=\2(ab+ba)$ is the Jordan product. Thus one obtains
an operator valued sesquilinear form
\begin{equation}\label{e2}
<a,b>=\phi(a\circ b^*)-\phi(a)\circ\phi(b)^*, a,b\in A.
\end{equation}
If we apply the Cauchy-Schwarz inequality to $\omega(<a,b>)$ for
all states $\omega$ of $B$ it was noticed in \cite{ER} that if
$\phi(a\circ a^*)= \phi(a)\circ\phi(a)^*$ then $<a,b>=0$ for all
$b\in A$. We call the set
$$
A_\phi =\{a\in A:\phi(a\circ a^*)= \phi(a)\circ\phi(a)^*\}
$$
the \textit{definite set} of $\phi$. It is a Jordan subalgebra of
$A$, and if $a\in A_\phi$ then $\phi(a\circ b)=\phi(a)\circ
\phi(b)$ for all $b\in A$.

In the present paper we shall develop the theory further.  We
first study positive unital normal, i.e. ultra weakly continuous,
maps $\phi\colon M \to M$, where $M$ is a von Neumann algebra. We
mainly study properties of the definite set $M_\phi$ and some of
its Jordan subalgebras of $M$ plus convergence properties of the
orbits $(\phi^n (a))$ for $a\in M$. We shall show that when there
exists a faithful family $\1 F$ of $\phi$-invariant normal states
there is a largest Jordan subalgebra $C_\phi$ of $M $ called the
multiplicative core of $M$, on which $\phi$ acts as a Jordan
automorphism. Furthermore if $a\in M$ then every weak limit point
of the orbit $(\phi^n (a))$ lies in $C_\phi$, and if $\rho(a\circ
b)=0$ for all $b\in C_\phi$, then $\phi^n (a)\rightarrow 0$
weakly.

Much of the above work was inspired by a theorem of Arveson,
\cite{A}. In the last section we study the $C^*$-algebra case and
the relation of our discussion with Arveson's work. Then $\phi
\colon A\to A$ is a positive unital map, and we assume the orbits
$(\phi^n (a))$ with $a\in A$ are norm relatively compact and that
there exists a faithful family $\1F$ of $\phi$-invariant states.
It is then shown that the multiplicative core $C_\phi$ of $\phi$
equals the set of main interest in \cite{A}, namely the norm
closure of the linear span of all eigenoperators $a\in A$ with
$\phi(a)=\lambda a, |\lambda|=1$, and that $
\lim_{n\rightarrow\infty} \|\phi^n (a)\|=0$ if and only if
$\rho(a\circ b)=0$ for all $b\in C_\phi$ and $\rho\in {\1F}$.

\section*{2.  Maps on von Neumann algebras}

Throughout this section $M$ denotes a von Neumann algebra,
$\phi\colon M\to M$ is a  positive normal unital map. $ M_\phi$
denotes the definite set of $\phi$ and $< , >$ the operator valued
sesquilinear form $<a,b>=\phi(a\circ b^*)-\phi(a)\circ\phi(b)^*,
a,b\in A$.
\begin{lem}\label{lem 2.1}
Let assumptions be as above, and suppose $(a_\alpha)$ is a bounded
net in $M$ which converges weakly to $a\in M$. If
$<a_\alpha,a_\alpha>\rightarrow 0$, then $a\in M_\phi$, and
$\phi(a\circ b)= \phi(a)\circ \phi(b)$ for all $b\in M.$
\end{lem}

\bp Let $\omega$ be a normal state on $M$. By the Cauchy-Schwarz
inequality, if $b,c\in M$ we have
$$
|\omega(<b,c>)|^2 \leq\omega(<b,b>)\omega(<c,c>).
$$
By assumption, if $a_\alpha$ and $a$ are as in the statement of
the lemma, and $b\in M$ then
\begin{eqnarray*}
|\omega(<a,b>)|^2 &=& \lim_\alpha |\ \omega(<a_\alpha,b>)|^2\\
&\leq &\lim_\alpha \omega(<a_\alpha,a_\alpha>)\omega(<b,b>) =0.
\end{eqnarray*}
Since this holds for all normal states $\omega$, $<a,b>=0$,
completing the proof.

In analogy with the definition of G-finite for automorphism groups
we introduce

\begin{defin}\label{def 2.2}
With $\phi$ as above we say $M$ is $\phi$-\textit{finite} if there
exists a faithful family $\1 F$ of $\phi$-invariant normal states
on the von Neumann algebra generated by the image $\phi(M)$.
\end{defin}

\begin{lem}\label{lem 2.3}
Assume $M$ is $\phi$-finite. Then for $a\in M$ we have
\enu{i}Every weak limit point of the orbit $(\phi^n (a))$ of $a$
belongs to $M_\phi$. \enu{ii}If $\rho(\phi^n (a)\circ b)=0$ for
all $b\in M_\phi, \rho \in \1F$, then $\phi^n (a)\rightarrow 0$
weakly.

\end{lem}
\bp If $\rho\in \1F$ denote by $\| . \|_\rho$ the seminorm
$\|x\|_\rho =\rho(x\circ x^*)^\2$. Then by the inequality (1)
\begin{eqnarray*}
\|\phi^{n+1}(a)\|_\rho ^2 &=&
\rho(\phi^{n+1}(a)\circ\phi^{n+1}(a)^*)\\
&\leq& \rho(\phi(\phi^n (a)\circ\phi^n (a)^*))\\
&=&\|\phi^n (a)\|_\rho ^2.
\end{eqnarray*}
Thus the sequence $\|\phi^n (a)\|_\rho ^2$ is decreasing, hence
$\|\phi^n (a)\|_\rho ^2 - \|\phi^{n+1}(a)\|_\rho ^2  \rightarrow
0$. We have
\begin{eqnarray*}
\rho(<\phi^n (a),\phi^n (a)>)&=&\rho(\phi(\phi^n (a))\circ\phi^n
(a)^*)-\phi(\phi^n (a))\circ\phi(\phi^n (a)^*))\\
&=&\rho(\phi^n (a)\circ\phi^n(a)^* -\phi^{n+1}
(a)\circ\phi^{n+1}(a)^*)\\
&=&\|\phi^n (a)\|_\rho ^2 - \|\phi^{n+1}(a)\|_\rho ^2  \rightarrow
0.
\end{eqnarray*}
Since this hold for all $\rho\in \1F$ and $\1F$ is faithful,
$<\phi^n (a),\phi^n (a)>\rightarrow 0$  weakly. By Lemma ~\ref{lem
2.1}, if $a_0$ is a weak limit point of $(\phi^n (a))$ then $a_0
\in M_\phi$, proving (i).

To show (ii) suppose $\rho(\phi^n (a)\circ b)=0$ for all $b\in
M_\phi, \rho\in \1F$. Let $a_0$ be a weak limit point of $(\phi^n
(a))$. Then $\rho(a_0 \circ b)=0$ for all $b\in M_\phi$, in
particular by part (i) $\rho(a_0\circ a_0)=0$. Since $\1F$ is
faithful on the von Neumann algebra generated by $\phi(M)$, $a_0
=0$. Thus 0 is the only weak limit point of $(\phi^n (a))$, so
$\phi^n (a)\rightarrow 0$ weakly.  The proof is complete.

It is not true in general that $\phi(M_\phi)\subseteq M_\phi$. We
therefore introduce the following auxiliary  concept.  If
$\phi\colon A\to A$ is positive unital with A a $C^*$-algebra,
then $A_{\Phi}=\{a\in A_{\phi} : {\phi}^k (a)\in A_{\phi}, k\in
\0N\}$.
\begin{lem}\label{lem 2.4}
Let $M$ be $\phi$-finite and $M_\Phi$ defined as above. Then
$M_\Phi$ is a weakly closed Jordan subalgebra of $M_\phi$ such
that $\phi(M_{\Phi})\subseteq M_{\Phi}$, and if $a\in M$ then every weak limit
point of $(\phi^n (a))$ belongs to $M_\Phi$. Furthermore, if
$\rho(\phi^n (a)\circ b)=0$ for all $b\in M_\Phi, \rho\in \1F$,
then $\phi^n (a)\rightarrow 0$ weakly.
\end{lem}

\bp Since $M$ is weakly closed and $\phi$ is weakly continuous on
bounded sets $M_\Phi$ is weakly closed.  Since $\phi$ and its
powers $\phi^k$ are Jordan homomorphisms on $M_phi$ it is
straightforward to show $M_\Phi$ is a Jordan subalgebra of $M$.
Furthermore it is clear from its definition that
$\phi(M_\Phi)\subseteq M_\Phi$.

If $a\in M$ and $a_0$ is a weak limit point of $(\phi^n (a))$, then
$a_0\in M_\phi$ by Lemma 3. Then $\phi(a_0)$ is a weak limit point
of $(\phi^{n+1}(a))$, hence belongs to $M_\phi$, again by Lemma 3.
Iterating we have $\phi^k (a_0)\in M_\phi$ for all $k\in\0N$. Thus
$a_0\in M_\Phi$. The last statement follows exactly as in Lemma 3.
The proof is complete.

It is not true that $\phi(M_\Phi)=M_\Phi$.  To remedy this problem
we introduce yet another Jordan subalgebra.

\begin{defin}\label{def 2.5}
Let $\phi\colon A\to A$ be positive unital with $A$ a
$C^*$-algebra. The \textit{multiplicative core} of $\phi$ is the
set
$$
C_\phi=\bigcap_{n=0}^{\infty}\phi^n (A_\Phi).
$$
\end{defin}
\begin{lem}\label{}
$C_{\phi}$ satisfies the following:
\enu{i} $C_{\phi}$ is a Jordan subalgebra of $A$.
\enu{ii} $\phi(C_{\phi})=C_{\phi}$.
\hfil\break
Suppose there exists a family $\1F$ of $\phi$-invariant states which is faithful on the $C^*$-algebra generated by $\phi(A)$. Then
we have
\enu{iii} The restriction of $\phi$ to $C_{\phi}$ is a Jordan automorphism.
\enu{iv} $C_{\phi}$ is the largest Jordan subalgebra of $A$ on which the restriction of $\phi$ is a Jordan automorphism.
\end{lem}
\bp 

As in Lemma 4 $C_\phi$ is clearly a Jordan subalgebra of $A$ such
that $\phi(C_\phi)\subseteq C_\phi$ and is weakly closed in the
von Neumann algebra case. Furthermore, since
$\phi(A_\Phi)\subseteq A_\Phi$, we have $\phi^n(A_\Phi)\subseteq
\phi^{n-1}(A_\Phi)$, so that the sequence $(\phi^n(A_\Phi))$ is
decreasing. Thus
$$
C_\phi =\bigcap_{n=0}^{\infty}\phi^{n+1}(A_\Phi)=\phi(C_\phi),
$$
so (i) and (ii) are proved.

We next show (iii), and let $\1F$ be as in the statement of the lemma. 
By (ii)  the restriction of $\phi$ to $C_\phi$ is a Jordan homomorphism 
of $C_{\phi}$ onto itself. In particular $\1F$ is faithful on $C_{\phi}$, so that  $\phi$ is
faithful on $C_{\phi}$, hence  is a Jordan automorphism of $C_\phi$, proving (iii).  

To show (iv)  let $B$ be a Jordan subalgebra of
$A$ such that $\phi|_B$ is a Jordan automorphism of $B$.  Then clearly
$B\subseteq A_\Phi$, and $\phi^n (B)=B$, so
that
$$
B=\bigcap_{n=0}^{\infty}\phi^n (B)
\subseteq\bigcap_{n=0}^{\infty}\phi^n (A_\Phi)=C_\phi.
$$
The proof is complete.
\hfil\break

We can now prove our main result.
\begin{thm}\label{Thm 2.5}
Let $M$ be $\phi$-finite, and $\1F$ a set of normal
$\phi$-invariant states which is faithful on the von Neumann
algebra generated by $\phi(M)$. Let $a\in M$. Then we have
\enu{i}Every weak limit point of $(\phi^n (a))$ lies in $C_\phi$.
\enu{ii} If $\rho(a\circ b)=0$ for all $b\in C_\phi, \rho\in \1F$,
then $\phi^n (a) \rightarrow 0$ weakly.
\end{thm}
\bp Ad(i). Let $a_0$ be a weak limit point of $(\phi^n (a))$.  By
Lemma 4 $a_0\in M_\Phi$. Choose a subnet $(\phi^{n_\alpha}(a))$
which converges weakly to $a_0$. Let $k\in \0N$, and let
$(\phi^{m_\beta} (a))$ be a subnet of $(\phi^{n_{\alpha} -k}(a))$
which converges weakly to $a_1\in M_\Phi$ (again using Lemma 4,
since $(\phi^{m_\beta} (a))$ will be a subnet of $(\phi^n (a))$).
Each $m_\beta$ is  of the form $n_{\alpha_j}-k$. The net
$(\phi^{n_{\alpha_j}}(a))$ converges to $a_0$, since it is a
subnet of the converging net $(\phi^{n_\alpha}(a))$. Thus we have
\begin{eqnarray*}
\phi^k (a_1)&=& \lim \phi^k (\phi^{m_\beta}(a))\\
&=& \lim \phi^{k+(n_{\alpha_{j}}-k)} (a)\\
&=& \lim \phi^{n_{\alpha_j}}(a)\\
&=& a_0.
\end{eqnarray*}
Thus $a_0\in \phi^k (M_\Phi)$ for all $k\in \0N$, hence $a_0\in
C_\phi$.

To show (ii) suppose $\rho(a\circ b)=0$ for all $\rho\in \1F,
k\in\0N$. Since $\phi^k (C_\phi)=C_\phi$ there exists $c\in
C_\phi$ such that $b=\phi^k (c)$. Thus
\begin{eqnarray*}
\rho(\phi^k (a)\circ b) &=& \rho(\phi^k (a)\circ \phi^k (c))\\
&=& \rho(\phi^k (a\circ b))\\
&=& \rho(a\circ b) = 0.
\end{eqnarray*}
By part (i) every weak limit point $a_0$ of $(\phi^n (a))$ lies in
$C_\phi$, so it follows by the above that $\rho(a_0 \circ b)=0$
for all $b\in C_\phi$. In particular $\rho(a_0 \circ a_0)=0$, so
by faithfulness of $\1F$, $a_0=0$, hence $\phi^n (a)\rightarrow 0$
weakly. The proof is complete.

One might believe that the converse of part (ii) in the above
theorem is true.  This is false. Indeed, let $M_0$ be a von
Neumann algebra with a faithful normal tracial state $\tau_0$. Let
$M_i=M_0, \tau_i =\tau_0, i\in \0Z$, and let $M=
\bigotimes_{-\infty}^{\infty} (M_i,\tau_i)$. Let $\phi$ be the
shift to the right. Then $C_\phi =M$. However, if $a= ...1\otimes
a_0 \otimes 1...\in M$ with $a_0\in M_0$, then
$\lim_{n\rightarrow\infty} \phi^n (a)= \tau_0 (a_0)1$, so if
$\tau_0(a_0)=0$, then the weak limit is 0. But $\tau(a\circ b)\neq
0$ for some $b\in M =C_\phi$.

If we assume convergence in the strong-* topology then the
converse holds, as we have

\begin{prop}\label{prop2.6}
Let $M$ be $\phi$-finite. Let $a\in M$ and suppose the sequence
$(\phi^n (a))$ converges in the strong-* topology. Then
$\rho(a\circ b)=0$ for all $b\in C_\phi, \rho\in \1F$ if and only
if $\phi^n (a)\rightarrow 0$ *-strongly.
\end{prop}
\bp If $\rho(a\circ b)=0$ for all $b\in C_\phi, \rho\in \1F$ then
$\phi^n (a)\rightarrow 0$ weakly by the theorem. Since the
sequence converges *-strongly the limit must be 0.

Conversely, if $\phi^n (a)\rightarrow 0$ *-strongly, then for all
$b\in C_\phi, \rho\in \1F$
$$
\rho(a\circ b)=\rho(\phi^n (a\circ b))=\rho(\phi^n (a)\circ\phi^n
(b))\rightarrow 0,
$$
since multiplication is *-strongly continuous on bounded sets. The
proof is complete.

We have not in general found a nice description of the complement
of $C_\phi$ in $M$, i.e. a subspace $D$ such that $M$ is a direct
sum of $C_\phi$ and $D$. In the finite case with a faithful normal
$\phi$-invariant trace this can be done.

\begin{prop}\label{prop 2.7}
Suppose $M$ has a faithful normal $\phi$-invariant tracial state.
Then there exists a faithful normal positive projection $P\colon
M\to C_\phi$ which commutes with $\phi$. Let $D=\{a-P(a): a\in
M\}$.  Then $M=C_\phi + D$ is a direct sum, and if $a\in D$ then
$\phi^n (a)\rightarrow 0$ weakly.
\end{prop}
\bp Since $M$ is finite the same construction as that of trace
invariant conditional expectations onto von Neumann subalgebras
yields the existence of a faithful trace invariant positive normal
projection $P \colon M\to C_\phi$, see \cite {HE}. Let $\tau$ be
the trace alluded to in the proposition.  Since $\tau$ is faithful
and $\phi$-invariant, $\phi$ has an adjoint map $\phi^* \colon
M\to M$ defined by $\tau(a \phi^*(b))=\tau(\phi(a)b)$ for $a,b\in
M$. Clearly $\phi^*$ is $\tau$-invariant, positive, unital, and
normal, and its extension $\bar\phi^*$ to an operator on
$L^2(M,\tau)$ is the usual adjoint of the extension $\bar\phi$ of
$\phi$. Since the restriction of $\bar\phi$ to the closure
$C_\phi^-$ of $C_\phi$ in $L^2(M,\tau)$ is an isometry of
$C_\phi^-$ onto itself, so is $\bar\phi^*$. It follows that $\phi
P=P\phi P=(P \phi^* P)^*=(\phi^* P)^* =P\phi$.

It is clear that $M=C_\phi + D$ is a direct sum. Suppose $a\in D$,
i.e. $P(a)=0$. Then $\tau(a\circ b)=0$ for all $b\in C_\phi$. If
we let $\1F =\{\tau|_{C_\phi}\circ P\}$ then, since $P$ commutes
with $\phi$, $\1F$ is a faithful family of normal $\phi$-invariant
states. By Theorem 6 $\phi^n (a)\rightarrow 0$ weakly, proving the
proposition.

\section*{3. Maps of C*-algebras}

Arveson \cite {A} proved the following result.

\begin{thm}\label{Thm.3.1}(Arveson)
Let $A$ be a $C^*$-algebra, $\phi\colon A\to A$ a completely
positive contraction such that the orbit $(\phi^n (a))$ is norm
relative compact for all $a\in A$. Then there exists a completely
positive projection $P\colon A\to A$ onto the norm closed linear
span $E_\phi$ of the eigenoperators $a\in A$ with
$\phi(a)=\lambda a$, with $|\lambda|=1$, and $\alpha
=\phi|_{E_\phi}$ is a complete isometry of $E_\phi$ onto itself.
We have
$$
\lim_{n\rightarrow\infty}\|\phi^n (a) - (\alpha\circ P)^n (a)\|=0,
$$
and $A$ is the direct sum of $E_\phi$ and the set $\{a\in A:\lim_n
\|\phi^n (a)\|=0\}$.
\end{thm}
We shall now show how our previous results yield a result which is
in a sense complementary to Arveson's theorem.

\begin {thm}\label{thm 3.2}
Let $A$ be a unital $C^*$-algebra and $\phi\colon A\to A$ a
positive unital map such that  the orbit $(\phi^n (a))$ is norm
relative compact for all $a\in A$. Let $C_\phi$ be the
multiplicative core for $\phi$ in $A$, and let $E_\phi$ denote the
set of eigenoperators $a\in A$ such that $\phi(a)=\lambda a$,
with $|\lambda|=1$. Assume there exists a set $\1F$ of
$\phi$-invariant states which is faithful on the $C^*$-algebra generated by $\phi(A)$. Then we have
 \enu{i} $E_\phi=C_\phi$ is a Jordan
subalgebra of $A$.
 \enu{ii} The restriction $\phi|_{E_\phi}$ is a Jordan automorphism
 of $E_\phi$.
 \enu{iii} Let $a\in A$. Then $\rho(a\circ b)=0$ for all $\rho\in\1F, b\in
 C_\phi$ if and only if $\lim_{n\rightarrow \infty}\|\phi^n
 (a)\|=0$.
 \end{thm}
\bp We first show (ii). If $\phi(a)=\lambda a$ then $\phi(a^*)=\bar\lambda a^*$, so
$E_\phi$ is self-adjoint. Furthermore by inequality (1)
$$
\phi(a\circ a^*)\geq \phi(a)\circ\phi(a^*)=\lambda a\circ \bar
\lambda a^*=a\circ a^*.
$$
Composing by $\rho\in\1F$ and using that $\1F$ is faithful on
$C^*(\phi(A))$ it follows that $\phi(a\circ
a^*)=\phi(a)\circ\phi(a^*)$, so $a\in A_\phi$, the definite set of
$\phi$. Since $a\in E_\phi$ is an eigenoperator, so is $a^2$,
hence $E_\phi$ is a Jordan subalgebra of $A_\phi$. Note that if
$\phi(a)=\lambda a$ then $\phi(\phi(a))=\phi(\lambda a)=\lambda  \phi(a))
$, so $\phi(a)\in E_\phi$. Thus $\phi\colon E_\phi \to E_\phi$.
If $a=\sum \mu_i a_i \in E_\phi$ where $\phi(a_i)=\lambda_i a_i$,
then $a=\sum \mu_i \bar\lambda{_i} \phi(a_i) \in \phi(E_\phi)$, so
by density of such $a'$s, $\phi(E_\phi)=E_\phi$. Thus by
faithfulness of $\1F$ the restriction $\phi|_{E_\phi}$ is a Jordan
automorphism, proving (ii).

It follows from Lemma 6 that $E_\phi\subseteq C_\phi$. To show the converse inclusion we
use that   the orbit $(\phi^n (a))_{n\in\0N}$ is norm relative
compact for all $a\in A$. By Lemma 6 the  restriction of $\phi$ to $C_\phi$ is a Jordan automorphism,
hence in particular an isometry. We assert that if $a\in C_\phi$
then the orbit $(\phi^n (a))_{n\in\0Z}$ is relative norm compact.
For this it is enough to show that the set $(\phi^{-n}
(a))_{n\in\0N}$ is relative norm compact, or equivalently that
each sequence $(\phi^{-n_k} (a))$ has a convergent subsequence. By
assumption $(\phi^{n_k} (a))$ has a convergent subsequence
$(\phi^{m_l} (a))$. Since this sequence is Cauchy, and
$$
\|\phi^{-n} (a) -\phi^{-m}
(a)\|=\|\phi^{n+m}(\phi^{-n}(a)-\phi^{-m} (a))\|=\|\phi^n (a)
-\phi^m (a)\|,
$$
it follows that $(\phi^{-m{_l}} (a))$ is Cauchy, and therefore
converges. Thus the set  $(\phi^{-n} (a))_{n\in\0N}$ is relative
norm compact, as is $(\phi^n (a))_{n\in\0Z}$. By a well known
result on almost periodic groups, see e.g. Lemma 2.8 in \cite{A},
$\phi|_{C_\phi}$ has pure point spectrum. Thus $C_\phi \subseteq
E_\phi$, proving (i).

It remains to show (iii). As in the proof of Lemma~\ref{lem 2.3}
we find that every norm limit point $a_0$ of $(\phi^n (a))$
belongs to $A_\phi$, and by the proof of Lemma ~\ref{lem 2.4}
$a_0\in A_\Phi =\{x\in A_\phi : \phi^k (x)\in A_\phi, k\in \0N$\}.
A straightforward modification of the proof of Theorem ~\ref{Thm
2.5}(i), replacing weak by norm, shows that $a_0\in C_\phi$. Let
$a\in A$ satisfy $\rho(a\circ b)=0$ for all $ b\in C_\phi, \rho\in
\1F$. Then by the proof of Theorem ~\ref{Thm 2.5}(ii), every norm
limit point of $(\phi^n (a))$ is 0. Thus there is a subsequence
$(\phi^{n_k}(a))$ of $(\phi^n (a))$ such that for all $\varepsilon
>0$ there is $k_0$ such that $\|\phi^{n_k}\|\leq \varepsilon$ when $k\geq k_0$.
But then $n> n_k$ for $k\geq k_0$ implies
$$
\|\phi^n (a)\|=\|\phi^{n-n_k}(\phi^{n_k}(a)\|\leq
\|\phi^{n_k}\|<\varepsilon.
$$
Thus $\|\phi^n (a)\|\rightarrow 0$.

Conversely, if  $\|\phi^n (a)\|\rightarrow 0$ then for $a\in A,
b\in C_\phi, \rho\in \1F$
$$
\rho(a\circ b)=\rho(\phi^n (a\circ b))=\rho(\phi^n (a)\circ \phi^n
(b))\rightarrow 0,
$$
 for  $n\rightarrow\infty$. Thus $\rho(a\circ
b)=0$, completing the proof of the theorem.

It was shown in \cite{ES} that if $A$ is a $C^*$-algebra, and
$P\colon A\to A$ is a faithful positive unital projection then the
image $P(A)$ is a Jordan subalgebra of $A$. The following
corollary proves more.

\begin{cor}\label{cor 3.3}
Let $A$ be a $C^*$-algebra and $P\colon A\to A$ a faithful
positive unital projection. Then $E_P =C_P =P(A)$, Hence $P(A)$ is
in particular a Jordan subalgebra of $A$.
\end{cor}
\bp Since $P^2 =P$ the orbit of each $a\in A$ is finite, so
compact. Since $P$ is faithful the set of states
$\1F=\{\omega|_{P(A)}\circ P\}$ with $ \omega$ a  state on $A$, is
a faithful family of $P$-invariant states. Thus by Theorem
~\ref{thm 3.2} we have $E_P = C_P$. Since $P$ is a projection the
only nonzero eigenvalue of $P$ is 1, and the corresponding eigen
operators are the elements in $P(A)$. Thus $E_P =P(A)$, proving
the corollary.

Department of Mathematics, University of Oslo, 0316 Oslo, Norway.

e-mail erlings@math.uio.no


\begin{thebibliography}{999}

\bibitem{A}
W. Arveson, \textsl{Asymptotic stability I: completely positive
maps}, Int.J.Math. 15(3) (2004), 289-312.

\bibitem{B}
B.M.Broise, letter to the author (1967).

\bibitem{C1}
M.-D. Choi, \textsl{Positive linear maps on $C^*$-algebras},
Thesis, University of Toronto (1972).

\bibitem{C2}
M.-D. Choi, \textsl{Positive linear maps of $C^*$-algebras},
Canad. J. Math. 24 (1972), 520-529.

\bibitem{ES}
E.Effros and E.St{\o}rmer, \textsl{Positive projections and Jordan
structure in operator algebras}, Math.Scand.45 (1979),127-138.

\bibitem{ER}
D.Evans and R.H{\o}egh-Krohn, \textsl{Spectral properties of
positive maps on $C^*$-algebras}, J.London Math.soc. 17(1978), 345-355.

\bibitem{HE}
U.Haagerup and E.St{\o}rmer, \textsl{Positive projections of von
Neumann algebras onto JW-algebras}, Reports on Math.Phys.36
(1995),317-330.

\bibitem{K1}
R.V.Kadison, \textsl{A general Schwarz inequality and algebraic
invariants for operator algebras}, Ann.Math.56(1952),494-503.

\bibitem{K2}
R.V.Kadison, \textsl{The trace in finite von Neumann algebras},
Proc.Amer.Math.Soc.12(1961),973-977.

\bibitem{S}
E.St{\o}rmer, \textsl{Positive linear maps of operator algebras},
Acta Math.110(1963),233-278.

\end{thebibliography}
\end{document}